\input{style/arxiv-general.cfg}
\documentclass[MSNbibl,number,citesort,seceqn,dvips]{arxbj}
\makeatletter
   \@ifpackageloaded{graphicx}{}{\usepackage{graphicx}}
\makeatother

\aid{0}
\volume{21}
\issue{2}
\pubyear{2015}
\firstpage{1002}
\lastpage{1013}
\doi{10.3150/14-BEJ595} % kopijuoti is 'New paper accepted'

\makeatletter
\newcommand{\rrvert}{\vert}
\newcommand{\llvert}{\vert}
\newtheorem{theorem}{Theorem}[section]
\newtheorem{lem}{Lemma}[section]
\newremark{rem}{Remark}
\newcommand{\E}{\mathbb{E}}
\newcommand{\Z}{\mathbb{Z}}
\newcommand{\R}{\mathbb{R}}
\newcommand{\nat}{\mathbb{N}}
\renewcommand{\Pr}{\mathbb{P}}
\newcommand{\ind}{\mathbf{1}}
\newcommand{\NB}{\operatorname{NB}}
\newcommand{\Ref}[1]{(\ref{#1})}
\newcommand{\Eq}{=}
\newcommand{\s}{\sigma}
\renewcommand{\th}{\theta}
\newcommand{\m}{\mu}
\newcommand{\D}{\Delta}
\renewcommand{\t}{\tau}
\newcommand{\cF}{\mathcal{F}}
\newcommand{\cL}{\mathcal{L}}
\newcommand{\Le}{\le}
\newcommand{\Def}{:=}
\newcommand{\ba}{\bar{a}}
\newcommand{\sji}{\sum_{j\ge1}}
\newcommand{\Po}{\operatorname{Po}}
\newcommand{\bR}{\overline{R}}
\renewcommand{\L}{\Lambda}
\newcommand{\hQ}{\widehat{Q}}
\renewcommand{\citet}{\cite}
\renewcommand{\citeyear}{\cite}
\newcommand{\keitimas}{\dvtx}
\makeatother

\begin{document}
\begin{frontmatter}

\title{Stein factors for negative binomial approximation in
Wasserstein distance}
\runtitle{Stein factors for negative binomial approximation in
Wasserstein distance}

\begin{aug}
%%%% inicialai - be tarpu
\author[A]{\inits{A.D.}\fnms{A.D.} \snm{Barbour}\thanksref{A}\ead[label=e1]{a.d.barbour@math.uzh.ch}},
\author[B]{\inits{H.L.}\fnms{H.L.} \snm{Gan}\thanksref{B,e2}\ead[label=e2,mark]{ganhl@ms.unimelb.edu.au}} \and
\author[B]{\inits{A.}\fnms{A.} \snm{Xia}\corref{}\thanksref{B,e3}\ead[label=e3,mark]{aihuaxia@unimelb.edu.au}}
%%\runauthor{} %% auto
\address[A]{Angewandte Mathematik,
Universit\"at Z\"urich,
Winterthurerstr.~190,
8057 Z\"urich, Switzerland.\\ \printead{e1}}
\address[B]{Department of Mathematics and Statistics, the University
of Melbourne,
Parkville, VIC 3010, Australia. \printead{e2,e3}}
\end{aug}

% HISTORY:
\received{\smonth{10} \syear{2013}}
\revised{\smonth{12} \syear{2013}}

% ABSTRACT
%
\begin{abstract}
The paper gives the bounds on the solutions to a Stein equation
for the negative binomial distribution that are needed for
approximation in terms of the
Wasserstein metric. The proofs are probabilistic, and follow the approach
introduced in Barbour and Xia (\textit{Bernoulli} \textbf{12} (2006) 943--954).
The bounds are used
to quantify the accuracy of negative binomial approximation to parasite counts
in hosts. Since the infectivity of a population can be expected to be
proportional
to its total parasite burden, the Wasserstein metric is the appropriate choice.
\end{abstract}

% KEYWORDS
% visi is mazosios raides ir pagal abecele
%
\begin{keyword}
\kwd{negative binomial approximation}
\kwd{Stein factors}
\kwd{Stein's method}
\kwd{Wasserstein distance}\vspace*{-2pt}
\end{keyword}

\end{frontmatter}

%s1 #&#
\section{Introduction}

The negative binomial distribution is widely used in biology to model
the counts of
individuals in populations, since such counts are frequently
overdispersed, making
the Poisson distribution an unsuitable choice. Indeed, the main
advantage of the negative binomial family over the Poisson family is
the extra flexibility
in fitting that results because the negative binomial family has a
second parameter.
However, for the distribution of parasites among hosts,
there are plausible mechanistic models \cite{K48,K52} that predict
a negative binomial distribution,
and it is of interest to know whether a member of the negative binomial
family would
still give a reasonable approximation, if the detailed assumptions of
such a model
were relaxed. One of the quantities of primary interest is then the
total rate of
output of infective stages, which can be expected to be closely related
to the
total number of parasites in the population \cite{KM}. Thus
the approximation needs to be
good when measured by a distance that limits the differences in expectation
of (not necessarily bounded) Lipschitz functionals, which makes the
Wasserstein metric
a natural choice. In this paper, we make negative binomial approximation
using Stein's method a practical proposition, by giving bounds on the solutions
of an appropriate Stein equation that correspond to Lipschitz test functions.

The negative binomial distribution $\NB(r,p)$ has probabilities given by
\[
\NB(r,p)\{k\} \Eq\frac{\Gamma(r+k)}{\Gamma(r)k!}(1- p)^r p^k,\qquad k
\in\Z_+:=\{ 0,1,\ldots\}; r>0, 0<p<1.
\]
One can check directly that $W \sim\NB(r,p)$ if and only if
\[
\E \bigl[ p(r+W)g(W+1) - Wg(W) \bigr] \Eq0
\]
for a sufficiently rich class of functions $g\keitimas\nat\to\mathbb
{R}$. One such class consists
of the solutions $g_f\keitimas\nat\to\R$ to the equations
%
%
%e1.1 #&#
\begin{equation}\label{steineq}
p(r+i)g_f(i+1) - ig_f(i) \Eq f(i) - \NB(r,p)\{f\},
\qquad f \in\cF_W,
\end{equation}
where $\cF_W := \{f\keitimas|f(x)-f(y)|\le|x-y|, \forall x,y\in\Z
_+\}
$ denotes the class of
Lipschitz functions on $\Z_+$, and $\NB(r,p)\{f\} := \E f(Z)$ for $Z
\sim\NB(r,p)$.
Then, for any random variable $W$ on $\Z_+$,
%
%
%e1.2 #&#
\begin{equation}
\label{Stein-expectation} \E f(W) - \NB(r,p)\{f\} \Eq\E \bigl[ p(r+W)g_f(W+1) -
Wg_f(W) \bigr],
\end{equation}
and, if we can bound the right-hand side of the above equation
uniformly for $f\in\cF_W$,
then we have a uniform bound for the left
hand side as well; but this corresponds precisely to a bound on the
Wasserstein distance
between $\cL(W)$ and $\NB(r,p)$.

In order to control the right-hand side of \Ref{Stein-expectation}, it
is typically necessary
to have bounds on the quantities
\[
G_1 \Eq\sup_{f \in\cF_W}\sup_{w \in\nat}g_f(w);
\qquad G_2 \Eq\sup_{f \in\cF_W}\sup_{w \in\nat}
\bigl\llvert g_f(w+1) - g_f(w) \bigr\rrvert .
\]
This note establishes the following result:

%
%
%th1.1 #&#
\begin{theorem}\label{mainthm}
For any $r > 0$ and $0 < p < 1$,
%
%
%e1.3 #&#
%e1.4 #&#
\begin{eqnarray}
\label{modeg}G_1 &\Eq&\frac{1}{1-p},
\\
\label{thmdeltag}G_2 &\Le&\min \biggl\{\frac{2}{1-p},
\frac
{1+p}{(1-p)^2}, \sqrt{\frac
{r_0}{rp(1-p)^3}} \biggr\},
\end{eqnarray}
where $r_0$ is the solution in $r > 1/2$ of the equation $\Gamma
(r-\frac{1}2)/\Gamma(r) = 3\sqrt{2\mathrm{e}}/8$,
and satisfies $\sqrt{r_0} \le3/2$.
\end{theorem}

The proof is given in Section~\ref{proofofmainthm}. In Section~\ref
{3}, we apply Theorem~\ref{mainthm} to approximating the distribution
of parasites in hosts.

%s2 #&#
\section{The proof of Theorem \texorpdfstring{\protect\ref{mainthm}}{1.1}}\label{proofofmainthm}

Setting $g_f(i) = h_f(i) - h_f(i-1)$, equation \Ref{steineq} becomes:
\[
f(i) - \NB(r,p)\{f\} \Eq p(r+i) \bigl(h_f(i+1) -
h_f(i)\bigr) - i\bigl(h_f(i) - h_f(i-1)
\bigr),
\]
where the right-hand side is the generator of an
immigration--birth--death process with constant immigration
rate $rp$, and \textit{per capita} birth and death rates $p$ and $1$,
respectively. More generally,
we let $Z_i := Z_i^{[a,b]}$ denote an immigration--birth--death process
with immigration rate $a$
and with \textit{per capita} birth and death rates $b$ and $1$,
respectively, having $Z_i(0) = i$.
We write $Y_i^{[b]}$ for $Z^{[0,b]}$.
% and $Z_i^{(r,p)}$ for
%$Z_i^{[r(1-p),1-p]}$.
% let $\{Y_i(t), t\ge0\}$ denote a similar process, but with
%immigration rate zero, having $Y_i(0)=i$.
From \citet{BX01},
%
%
%e2.1 #&#
\begin{equation}
\label{solutionh} h_f(i) \Eq-\int_0^\infty
\bigl[ \E f\bigl(Z_i^{[rp,p]}(t)\bigr) - \NB(r,p)\{ f\} \bigr]
\,\mathrm{d}t.
\end{equation}

We make use of the following two lemmas, proved in \citet{K48}, who
attributes the
first to Palm. We write
%
%
%e2.2 #&#
\begin{equation}
\label{ADB-Lambda-theta-def} \L_t(b) \Def\mathrm{e}^{-(1-b)t} \quad\mbox{and}
\quad\th_t(b) \Def1 - (1-b)/\bigl(1 - b\L_t(b)\bigr).
\end{equation}

%
%
%le2.1 #&#
\begin{lem}\label{Y1}
$Y_1^{[b]}(t)$ has a modified geometric distribution: for $0 < b \neq1$,
\[
\Pr\bigl[Y_1^{[b]}(t) = 0\bigr] \Eq b^{-1}
\th_t;\qquad\Pr\bigl[Y_1^{[b]}(t) = k\bigr] \Eq
\L _t(1-\th_t)^2 \th_t^{k-1},
\qquad k\ge1,
\]
where $\L_t = \L_t(b)$ and $\th_t = \th_t(b)$. In particular, the
first two moments are given by
\[
\E Y_1^{[b]}(t) \Eq\L_t;\qquad \E\bigl
\{Y_1^{[b]}(t)\bigr\}^2 \Eq\frac{\L_t(1+b-2b\L_t)}{1-b}.
\]
If $b=1$, the limiting formulae as $b\to1$ hold true; for instance,
$\th_t(1) = t/(1+t)$ and
$\E\{Y_1^{[1]}(t)\}^2 = 1 + 2t$.
\end{lem}

%
%
%le2.2 #&#
\begin{lem}\label{NbZ0}
$Z_0^{[a,b]}(t)$ has the negative binomial distribution $\NB(a/b,\th_t)$.
%, with parameters $r$ and $p_t^*=\frac{p}{1-(1-p)e^{-pt}}$.
\end{lem}

\begin{pf*}{Proof of Theorem~\ref{mainthm}}
As $g_f(i) = h_f(i) -
h_f(i-1)$, it follows
from \Ref{solutionh} that
\[
g_f(i) \Eq-\int_0^\infty\E \bigl[ f
\bigl(Z_i(t)\bigr) - f\bigl(Z_{i-1}(t)\bigr) \bigr] \,
\mathrm{d}t,
\]
where, throughout the proof, we write $Z_j$ for $Z_j^{[rp,p]}$.
We now couple $Z_{i-1}$ and $Z_i$ by setting
\[
Z_i(t) \Eq Z_{i-1}(t) + Y_1(t),
\]
where $Y_1 \stackrel{d}{=} Y_1^{[p]}$, and $Z_{i-1}(t)$ and $Y_1(t)$
are independent. Then $g_f(i)$ can
be expressed as
\[
g_f(i) \Eq-\int_0^\infty\E \bigl[ f
\bigl(Z_{i-1}(t) + Y_1(t)\bigr) - f\bigl(Z_{i-1}(t)
\bigr) \bigr] \,\mathrm{d}t.
\]
Now, because $f \in\cF_W$, it follows that
\[
\bigl |g_f(i)\bigr | \Le\int_0^\infty\E
Y_1(t) \,\mathrm{d}t \Eq\int_0^\infty
\L_t(p) \,\mathrm{d}t \Eq\frac{1}{1-p},
\]
using Lemma~\ref{Y1} for the first equality, and this maximal value
for $|g_f|$ is attained by taking
$f(x)=-x$. This completes the proof of \Ref{modeg}, and also yields
the bound $2/(1-p)$ in $G_2$.

To prove the remainder of \Ref{thmdeltag}, we first observe that the
function that
maximizes $\Delta g_f(i)$ is $f_i(j) = -|j-i|$. This follows by using
the same argument as in
% The argument closely follows the one in
\citeyear{BX02}, proof of (1.4). In the rest of the proof, we
write $f=f_{i}$. Using the couplings
\[
Z_{i+1}(t) \Eq Z_i(t) + Y_1(t);\qquad
Z_i(t) \Eq Z_{i-1}(t) + Y_1'(t),
\]
where $Y_1,Y_1' \stackrel{d}{=} Y_1^{[p]}$ and the processes
$Z_{i-1}$, $Y_1$ and $Y_1'$ are independent, we obtain
\begin{eqnarray*}
\Delta g_f(i) &=&\int_0^\infty\E
\bigl[ f\bigl(Z_{i+1}(t)\bigr) - f\bigl(Z_i(t)\bigr) + f
\bigl(Z_{i-1}(t)\bigr) - f\bigl(Z_i(t)\bigr) \bigr] \,
\mathrm{d}t
\\
&\leq&\int_0^\infty\E\bigl\{ \bigl | Z_{i-1}(t) +
Y_1(t) + Y_1'(t) - i\bigr | -
\bigl |Z_{i-1}(t) + Y_1(t)-i\bigr |
\\
&&\phantom{\int_0^\infty\E\bigl\{} {}
-\bigl |Z_{i-1}(t) + Y_1'(t)-i\bigr | +
\bigl |Z_{i-1}(t)-i\bigr | \bigr\} \,\mathrm{d}t
\\
&\le&\int_0^\infty\E\bigl[2\min{
\bigl(}Y_1(t),Y_1'(t){\bigr)}
\ind_{\{
Z_{i-1}(t)<i<Z_{i-1}(t)+Y_1(t)+Y_1'(t)\}}\bigr] \,\mathrm{d}t,
\end{eqnarray*}
where the last inequality is because the quantity in the braces equals
$0$ if $Z_{i-1}(t)\ge i$ or $Z_{i-1}(t)+Y_1(t)+Y_1'(t)\le i$; it is
bounded by $2Y_1'(t)$ if one applies the triangle inequality to
$| Z_{i-1}(t) + Y_1(t) + Y_1'(t) - i| - |Z_{i-1}(t) + Y_1(t)-i|$ and
$-|Z_{i-1}(t) + Y_1'(t)-i| + |Z_{i-1}(t)-i|$, and hence it is also
bounded by $2Y_1(t)$ if one swaps $Y_1(t)$ and $Y_1'(t)$. This implies that
%
%
%e2.3 #&#
\begin{eqnarray}
\label{ADB-cross} \Delta g_f(i) %&\le\int_0^\infty\E\Big{[}2\min\{Y_1(t),Y_1'(t)\}
&\Le&\int
_0^\infty\sum_{i_1,i_2,j} 2
\min(i_1,i_2)\ind_{\{
i+1-i_1-i_2\le j\le i-1\}}
\nonumber
\\
&&\phantom{\int_0^\infty\sum
_{i_1,i_2,j}} {}\times\Pr \bigl(Y_1(t)=i_1
\bigr)\Pr\bigl(Y_1'(t)=i_2\bigr)\Pr
\bigl(Z_{i-1}(t)=j\bigr) \,\mathrm{d}t
\nonumber
\\
&\Le&\int_0^\infty\max_{j}\Pr
\bigl(Z_{i-1}(t)=j\bigr)
\\
&&\phantom{\int_0^\infty} {} \times\sum
_{i_1,i_2} 2\min(i_1,i_2)
(i_1+i_2-1) \Pr\bigl(Y_1(t)=i_1
\bigr)\Pr\bigl(Y_1'(t)=i_2\bigr) \,
\mathrm{d}t
\nonumber
\\
&\le& \int_0^\infty\max_{j}
\Pr\bigl(Z_{i-1}(t)=j\bigr)\E \bigl[\bigl(Y_1(t)+Y_1'(t)
\bigr) \bigl(Y_1(t)+Y_1'(t)-1\bigr)\bigr]
\,\mathrm{d}t.
\nonumber
\end{eqnarray}
To bound $\Pr[Z_{i-1}(t)=j]$, we decompose $Z_{i-1}(t)$ into a sum of
two independent components
\[
Z_{i-1}(t) \stackrel{d} {=} Z_0(t) + Y_{i-1}(t),
\]
where $Y_{i-1}\stackrel{d}{=}Y_{i-1}^{[p]}$, as defined earlier. From
this it follows,
using Lemma~\ref{NbZ0}, that
%
%
%e2.4 #&#
\begin{equation}
\label{xiaadd01} \max_j \Pr\bigl[Z_{i-1}(t) = j\bigr]
\Le\max_k \Pr\bigl[Z_0(t)=k\bigr] \Eq P(r,\th
_t),
\end{equation}
where $\th_t = \th_t(p)$ and $P(r,q) := \max_k \NB(r,q)\{k\}$.
In \citet{Phillips}, the representation of $\NB(r,q)$ as a
$\Gamma(r,(1-q)/q)$ mixed Poisson distribution,
where $\Gamma(r,\lambda)$ denotes the Gamma distribution with shape
parameter $r$ and scale parameter $1/\lambda$,
is exploited to bound $P(r,q)$. Using the bound $\max_k \Po(\lambda
)\{k\} \le1/\sqrt{2\mathrm{e}\lambda}$
from \citet{BJ89}, he shows that, if $r > 1/2$, then
\[
P(r,q) \Le\sqrt{\frac{1-q}{2\mathrm{e}rq}} K_r,
\]
where $K_r := \sqrt r \Gamma(r-\frac{1}2)/\Gamma(r)$ is decreasing in
$r > 1/2$.
Hence, since
\[
\frac{1-\th_t}{\th_t} \Eq\frac{1-p}{p(1-\L_t)},
\]
we have, for $\th_t = \th_t(p)$ and $\L_t = \L_t(p)$,
\[
P(r,\th_t) \leq %
\cases{ \sqrt{\displaystyle
\frac{1-p}{2\mathrm{e}rp(1-\L_t)}} K_r, &\quad$\mbox{if }r > 1/2$;\vspace*{4pt}
\cr
1,&\quad$
\mbox{if } r \leq1/2$. } %
\]

For the third element in the bound \Ref{thmdeltag}, we assume that $r
> 1/2$, and use \Ref{ADB-cross} to give
\begin{eqnarray*}
\Delta g_f(i) &\Le&\int_0^\infty\E {
\bigl[} \bigl(Y_1(t) + Y_1'(t) \bigr)
\bigl(Y_1(t) + Y_1'(t)-1\bigr) {\bigr]}
\sqrt{\frac{1-p}{2r\mathrm{e}p}} K_r \frac{1}{\sqrt{1-\L_t}} \, \mathrm{d}t
\\
& \Eq&\sqrt{\frac{1-p}{2r\mathrm{e}p}} K_r \cdot 2\int_0^\infty
\frac{\L_t((1-3p)\L_t +2p)}{(1-p) \sqrt{1-\L_t}} \, \mathrm{d}t,
\end{eqnarray*}
using the moments given in Lemma~\ref{Y1}.
Direct computations now give
\[
\int_0^\infty\frac{\L_t}{\sqrt{1-\L_t}} \,\mathrm{d}t \Eq
\frac
{2}{1-p},\qquad \int_0^\infty
\frac{\L_t^2}{\sqrt{1-\L_t}} \,\mathrm{d}t \Eq \frac{4}{3(1-p)},
\]
leading to the result
%
%
%e2.5 #&#
\begin{equation}
\label{ADB-hard-bnd} \|\Delta g_f\| \Le\frac{8}{3} \sqrt{
\frac{1}{2r\mathrm
{e}p(1-p)^3}} K_r, \qquad r > 1/2. %\Le\frac3{\sqrt{2rp^3(1-p)}}.
\end{equation}
Note that, for any $p$, this is at least $16K_r/\{3(1-p)\sqrt{2\mathrm
{e}r}\}$,
which is smaller than $2/(1-p)$ whenever
$r > r_0$, for $r_0^{-1/2}K_{r_0} = 3\sqrt{2\mathrm{e}}/8$. Hence,
%
%
%e2.6 #&#
\begin{equation}
\label{ADB-r_0-bnd} \|\Delta g_f\| \Le\min \biggl\{\frac{2}{1-p},
\sqrt{\frac
{r_0}{rp(1-p)^3}} \biggr\},
\end{equation}
and computation gives $ \sqrt{r_0} \le1.427 < 3/2$.

Finally, for any $p$, $r$, we can simply bound $\max_j \Pr[Z_{i-1}(t) =
j]$ by $1$ in \Ref{ADB-cross}, giving
\begin{eqnarray*}
\|\Delta g_f\| &\Le&\int_0^\infty\E
\bigl[ \bigl(Y_1(t) + Y_1'(t) \bigr)
\bigl(Y_1(t) + Y_1'(t)-1\bigr)\bigr] \,
\mathrm{d}t
\\
&\Le&2\int_0^\infty\frac{\L_t((1-3p)\L_t +2p)}{1-p} \,\mathrm
{d}t \Eq\frac
{1+p}{(1-p)^2}.
\end{eqnarray*}
This bound is valid irrespective of the choices of $r>0$ and $0 < p < 1$.
% For $r\le1$, some improvement
% can be obtained by using $P(r,p_t^*) \le(p_t^*)^r$. For $r<1$ fixed,
%the bound is then of order
% $O(p^{r-2})$ as $p \to0$.
\end{pf*}

%%
%%
%In the case where $r \leq\frac{1}{2}$, if we use $p^*_t$ instead of
%$1$, a slightly sharper bound
%can be obtained as
%%
%(3p-2)(p^{r-2}-1)(1-r)\right]}{(1-p)^2(1-r)(2-r)}.
%%
%%
%}

%
%
\begin{rem*}
Note that the bounds in Theorem~\ref{mainthm} correspond exactly to
the bounds derived in \citet{BX02},
in the limit when $rp \to\lambda$ and $p \to0$, giving the Poisson case.
\end{rem*}

%s3 #&#
\section{An application to a parasite model}\label{3}

The model that we use to describe the development over time of the
number of parasites in
a host is based on the immigration--birth--death process $Z_0^{[a,b]}$ of
the previous section,
% with $r = a/(1-p)$,
with $a$ the rate of ingestion of parasites and $b$ their \textit{per
capita} birth rate. This model
would imply exactly negative-binomially distributed parasite numbers in
any age class.
However, since in reality $a$ can be expected to be variable, both
between individuals
and over time, we replace it by a function $a_t$, and investigate how
much this
influences the distribution of the number $W$ of parasites at some
fixed age $T$.
We fix any $\ba> 0$, to be thought of as a typical parasite ingestion
rate, and define
\begin{eqnarray*}
A_t &:=& \int_0^t
(a_{T-s}-\ba)\mathrm{e}^{-(1-b)s} \,\mathrm {d}s;\qquad
A_T^* := \sup_{0\le
t\le T} |A_t|;
\\
% \th_T &:=& \frac{b(1-e^{-(1-b)T})}{1-be^{-(1-b)T}};
R_T &:=& \frac{1-b}{b(1-\mathrm{e}^{-(1-b)T})} \int_0^T
a_{T-s}\mathrm {e}^{-(1-b)s} \,\mathrm{d}s,
\end{eqnarray*}
also setting $\th_T = \th_T(b)$ and $R_a^* := \ba/b$.
$A_t$ is a measure of the amount by which the cumulative exposure at
time $t$
under an ingestion rate of $a_s$, $0\le s\le t$, differs from that with constant
ingestion rate $\ba$, allowing for the evolution of the parasites
between ingestion
and time $T$.
Thus, both $|A_t|$ and $A_T^*$ reflect how closely the choice of $\ba$
corresponds to the actual
ingestion rate. If $R_T = R_a^*$, then $A_T=0$.

%
%
%th3.1 #&#
\begin{theorem} \label{example2-thm}
Under the above circumstances,
we have
\begin{eqnarray*}
% \e_1 &=& \frac{(1-\th_T)|A_T|}{b(1-b)(1-e^{-(1-b)T})}; \non\\
&&d_W\bigl(\NB\bigl(R_a^*,
\theta_T\bigr),W\bigr)
\\
&&\quad\le|A_T| + {16}\th_T A_T^*\bigl(1+
\ln\bigl\{1/(1-\th_T)\bigr\}\bigr) \min \biggl\{\frac{2}{1-\th_T},
\frac{3}{2\sqrt{R_a^*\theta
_T(1-\theta
_T)^3}} \biggr\}. % \label{example2-01}
\end{eqnarray*}
%
%$K$ can be taken to be $16$.
% In particular, if $\ba$ is such that $R_a = R_T$, $A_T=0$, and the
%bound simplifies somewhat.
%$b=0$, we have $W\stackrel{d}{=}\Pn(\lambda)$ with $\lambda=
\end{theorem}

\begin{rem*}
If $\{a_s\}$ includes a random component, $|A_T|$ and $A_T^*$ should be
replaced by their expectations
in the bound given in the theorem.
\end{rem*}

\begin{pf*}{Proof of Theorem~\ref{example2-thm}}
We define $N :=\{N_s, 0\le s\le T\}$ to be a Poisson process with mean function
$\E N_t = \int_0^t a_u \,\mathrm{d}u$. Given that the points of $N$
in $[0,T]$
are $\t_1 < \t_2 < \cdots$\,,
we sample values $(X_j, j\ge1)$ independently from the distributions
$\cL(Y_1^{[b]}(T-\t_j))$,
% with $p$ replaced by $(1-b)$,
and let $\Xi$ be the point process with $\Xi\{(0,s)\} := \sum_{j\keitimas\t_j < s} X_j$.
Then $W\stackrel{d}{=}\int_0^T \Xi(\mathrm{d}s)$.
For each $f\in\cF_W$, let $g:=g_f$ be a solution to the Stein
equation \Ref{steineq}
with $p = \th_T$ and $r = R_a^*$.
Since
\[
\cL\bigl(W |N\{s\} = 1\bigr) \Eq\cL\bigl(W + Y_1^*(T-s)\bigr),
\]
where $Y_1^* \stackrel{d}{=} Y_1^{[b]}$ is independent of $W$, we have
\begin{eqnarray*}
\E Wg(W)&=& \E \biggl\{\int_0^T g \bigl(\Xi
\bigl\{[0,T]\bigr\} \bigr)\Xi (\mathrm{d}s) \biggr\}
% &=&\E\int_0^T\eta^\ast_sg\left(\int_0^T\eta_tdN_t+\eta^\ast_s
 \Eq \sum_{j\ge1}j\E g(W+j)
\int_0^T\Pr\bigl[Y_1^{[b]}(T-s)=j
\bigr]a_s \,\mathrm{d}s
\\
&\Eq&\sum_{j\ge1}\E g(W+j)jC_j^T,%\label{example2-02}
\end{eqnarray*}
where $C_j^T:=\int_0^T \Pr[Y_1^{[b]}(T-s)=j]a_s \,\mathrm{d}s$.
Hence, for any $r$,
%
%
%e3.1 #&#
\begin{eqnarray}
\label{example2-03} &&\E \bigl(\theta_T(r+W)g(W+1)-Wg(W) \bigr)
\nonumber
\\
&&\quad= r\theta_T \E g(W+1)+\theta_T\sum
_{j\ge1}\E g(W+j+1)jC_j^T-\sum
_{j\ge1}\E g(W+j)jC_j^T
\\
&&\quad= \bigl(r\theta_T -C_1^T\bigr)\E
g(W+1)+\sum_{j\ge2}\E g(W+j) \bigl(\theta
_T(j-1)C_{j-1}^T-jC_j^T
\bigr).
\nonumber
\end{eqnarray}
Using Lemma~\ref{Y1}, we can verify that
\begin{eqnarray*}
\sum_{j\ge1}jC_j^T &=& \int
_0^T a_s\sum
_{j\ge1}j\Pr \bigl[Y_1^{[b]}(T-s)=j\bigr]
\,\mathrm{d}s \Eq\int_0^T a_s\E
Y_1^{[b]}(T-s) \,\mathrm{d}s
\\
&=& \int_0^T a_s
\mathrm{e}^{-(1-b)(T-s)} \,\mathrm{d}s,
\end{eqnarray*}
which in turn implies that
\[
-\sum_{j\ge2}\bigl(\theta_T(j-1)C_{j-1}^T-jC_j^T
\bigr) \Eq(1-\th_T)\sji jC_j^T -
C_1^T \Eq r \theta_T -C_1^T,
\]
if $r = R_T$.
%We first consider approximation by the distribution $\NB(R_T,1-
% R_T := \th_T^{-1}(1-\th_T)\int_0^T e^{-(1-b)(T-s)}a_s ds$.
Thus, it follows from \Ref{example2-03} that
%
%
%e3.2 #&#
\begin{eqnarray}
\label{example2-05} &&\E \bigl(\theta_T(R_T+W)g(W+1)-Wg(W)
\bigr)
\nonumber
\\[-8pt]
\\[-8pt]
&& \quad\Eq\sum_{j\ge2}\bigl(\E g(W+j)-\E g(W+1)\bigr)
\bigl(\theta _T(j-1)C_{j-1}^T-jC_j^T
\bigr).
\nonumber
\end{eqnarray}

On the other hand, Lemma~\ref{Y1} shows that
%
%
%e3.3 #&#
\begin{equation}
\label{Y_1-probs} \Pr\bigl[Y_1^{[b]}(t) = j\bigr] \Eq(1-
\th_t)^2 \mathrm{e}^{-(1-b)t} \th_t^{j-1}.
\end{equation}
Hence, defining $\bar C_j^T := \bar a\int_0^T\Pr[Y_1^{[b]}(T-s)=j]
\,\mathrm{d}s$, it follows that
$ \bar C_j^T \Eq\bar a\theta_T^j/(jb)$, $j\ge1$, which in turn gives
%
%
%e3.4 #&#
\begin{equation}
\label{example2-08} (j-1)\theta_T\bar C_{j-1}^T-j\bar
C_j^T \Eq0.
\end{equation}
Combining \Ref{example2-05} and \Ref{example2-08} and using Lemma~\ref{Y1} yields
%
%
%e3.5 #&#
\begin{eqnarray}
\label{ADB-pause} && \bigl\llvert \E \bigl(\theta_T(R_T+W)g(W+1)-Wg(W)
\bigr)\bigr\rrvert
\nonumber
\\
&&\quad\le \|\Delta g\|\sum_{j\ge2}(j-1) \bigl |
\theta_T(j-1) \bigl(C_{j-1}^T-\bar
C_{j-1}^T \bigr)-j \bigl(C_j^T-\bar
C_j^T \bigr)\bigr  |
\nonumber
\\[-8pt]
\\[-8pt]
&&\quad= \|\Delta g\|\sum_{j\ge2}(j-1)
\nonumber
\\
&&\phantom{\quad= \|\Delta g\|\sum_{j\ge2}}{} \times \Biggl|\int_0^T
(a_{T-s}-\bar a) \bigl(\theta_T(j-1)\Pr\bigl[Y_1^{[b]}(s)=j-1
\bigr] - j\Pr\bigl[Y_1^{[b]}(s)=j\bigr]\bigr) \, \mathrm{d}s
\Biggr|,
\nonumber
\end{eqnarray}
which, with \Ref{Y_1-probs}, allows concrete estimates to be undertaken.

The simplest and most direct
strategy is to impose bounds on $|a_{T-s}-\ba|$. However, this may not
lead to practically
useful results. For instance, animals may sleep at night and graze
during the day, so that $a_s$
can have substantial variation, but over time scales typically much
faster than the life history of
the parasite. Instead, we prefer to formulate bounds expressed in terms
of differences between cumulative exposure, which may more reasonably
be expected to be small.
For this reason, we write the quantity within the moduli in \Ref
{ADB-pause} as
\[
\int_0^T (a_{T-t}-\ba) \bigl(
\th_T(j-1)-j\th_t\bigr)\th_t^{j-2}(1-
\th_t)^2 \mathrm{e}^{-(1-b)t} \,\mathrm{d}t,
\]
write $f_j(\th) := (\th_T(j-1)-j\th)\th^{j-2}(1-\th)^2$ and
integrate by parts, giving
%
%
%e3.6 #&#
\begin{equation}
\label{ADB-integrated-j-bnd} A_T f_j(\th_T) - \int
_0^T A_t f_j'(
\th_t)\frac{\mathrm{d}\th
_t}{\mathrm{d}t} \,\mathrm{d}t,
% \Eq A_T f_j(\th_T) - \int_0^{\th_T} f_j'(\th) d\th,
\end{equation}
where $A_t := \int_0^t (a_{T-s}-\ba)\mathrm{e}^{-(1-b)s} \,\mathrm{d}s$.
Now the first term in \Ref{ADB-integrated-j-bnd} can easily be
bounded, because
$|f_j(\th_T)| = (1-\th_T)^2 \th_T^{j-1}$.
For the second, we use the bound
%
%
%e3.7 #&#
\begin{equation}
\label{ADB-second-bnd} \biggl\llvert \int_0^T
A_t f_j'(\th_t)
\frac{\mathrm{d}\th_t}{\mathrm
{d}t} \,\mathrm{d}t \biggr\rrvert \Le A_T^* \int
_0^{\th_T} \bigl |f_j'(\th)\bigr | \,
\mathrm{d}\th.
\end{equation}
Observe that
%
%
%e3.8 #&#
\begin{equation}
\label{f_j-dash} f_j'(\th) \Eq\th^{j-3}\bigl\{(1-
\th_T) + (\th_T - \th)\bigr\} Q_j\bigl(\th
_T,j(\th_T-\th),j(1-\th_T)\bigr),
\end{equation}
where $Q_j(x,y,z)$ is a homogeneous multinomial of degree $2$ in its
arguments and has coefficients
that are uniformly bounded in $j\ge2$, with the coefficient of $z^2$
being zero. Hence,
$|f_j'(\th)|$ can be bounded above by replacing $Q_j$ by $\hQ_j$ in
\Ref{f_j-dash}, where $\hQ_j$ is obtained
from $Q_j$ by taking the absolute values of its coefficients.
Integrating any of the terms from $0$ to $\th_T$ gives a bounded multiple
of either $j^{-2}\th_T^{j+1}$, $j^{-1}\th_T^j(1-\th_T)$ or $\th
_T^{j-1}(1-\th_t)^2$
to go into \Ref{ADB-second-bnd}, and multiplying each of these by
$(j-1)$ and adding
over $j\ge2$, as required by \Ref{ADB-pause}, gives a multiple of
$\th_T^2\log\{1/(1-\th_T)\}$, $\th_T^2$ or $\th_T$, respectively.
Hence, it follows
that
%
%
%e3.9 #&#
\begin{eqnarray}
\label{ADB-bound-with-delta} &&\bigl\llvert \E \bigl(\theta_T(R_T+W)g(W+1)-Wg(W)
\bigr)\bigr\rrvert
\nonumber
\\
&&\quad\Le\|\Delta g\| \th_T \bigl\{|A_T| +
A_T^*\bigl(K_1 + K_2 \log\bigl\{1/(1-
\th_T)\bigr\}\bigr)\bigr\}
\\
&&\quad\Le K\|\Delta g\| \th_T A_T^*\bigl(1 + \log\bigl
\{1/(1-\th_T)\bigr\}\bigr),
\nonumber
\end{eqnarray}
for suitable constants $K_1$, $K_2$ and $K$.
Careful computation in the \hyperref[appendix]{Appendix} shows that
$K_1 \le34/3$ and $K_2
\le16$, giving $K\le16$.

We now use \Ref{thmdeltag} of Theorem~\ref{mainthm} to bound $\|\D g\|
$ for all $g = g_f$, where $f \in\cF_W$ and
$g_f$ satisfies \Ref{steineq} with $r = R_a^*$ and $p = \th_T$; in
particular, this gives
\[
\|\Delta g\|\le\min \biggl\{\frac{2}{1-\th_T}, \frac{3}{2\sqrt {R_a^*\theta_T(1-\theta_T)^3}} \biggr\}.
\]
Therefore, it follows from \Ref{ADB-bound-with-delta} that
%
%
%e3.10 #&#
\begin{eqnarray}
\label{R_T-bnd} &&\bigl\llvert \E \bigl(\theta_T(R_T+W)g(W+1)-Wg(W)
\bigr)\bigr\rrvert
\nonumber
\\[-8pt]
\\[-8pt]
&&\quad\Le{16}\th_T A_T^*\bigl(1 + \ln\bigl\{1/(1-
\th_T)\bigr\}\bigr) \min \biggl\{ \frac{2}{1-\th_T},
\frac{3}{2\sqrt{R_a^*\theta_T(1-\theta_T)^3}} \biggr\}.
\nonumber
\end{eqnarray}
%
% completing the approximation by $\NB(R_T,1-\th_T)$.
But it is immediate from \Ref{Stein-expectation} and \Ref{modeg} that
\[
d_W\bigl(\NB(R_T,\th_T),\NB
\bigl(R_a^*,\th_T\bigr)\bigr) \le\frac{\th_T}{1-\th
_T}\bigl |R_T
- R_a^*\bigr | \Eq|A_T|,
\]
completing the proof of the theorem.
\end{pf*}

%
%
%re1 #&#
\begin{rem*}
Note also that, if $b\to0$ while $a$ is held fixed, then
$\theta_T \asymp b \to0$,
% and $R_a\theta_T$ %\to\m_T := \int_0^T
%e^{-(T-s)}a_s ds$,
% is bounded away from zero,
so that the upper bound in \Ref{R_T-bnd} approaches $0$. In this
limiting case, the number of parasites
has precisely a Poisson distribution, even for time varying $a$, with
mean $\m_T := \int_0^T \mathrm{e}^{-(T-s)}a_s \,\mathrm{d}s$.
\end{rem*}

Similar considerations can be applied to the distribution of the total
parasite burden
$W := \sum_{i=1}^nW^{(i)}$ among $n$
independent individuals, with their own functions $a^{(i)}$, $1\le i\le
n$, but all with the
same $b$. First, defining $\bR:= n^{-1}\sum_{i=1}^nR_T^{(i)}$, it follows
easily from \Ref{ADB-bound-with-delta} that
\[
d_W \bigl(\cL(W),\NB(n\bR,\th_T) \bigr) \Le\sup
_{f\in\cF_W}\|\Delta g_f\| {16}\th_T \sum
_{i=1}^n \bigl(A_T^*
\bigr)^{(i)}\bigl(1 + \ln\bigl\{1/(1-\th_T)\bigr\}\bigr),
\]
where $g = g_f$ satisfies \Ref{steineq}, with $r = n\bR$ and $p = \th
_T$, also because
$A_T^{(i)}= 0$ when approximating by $\NB(R_T^{(i)},\th_T)$. Hence,
for example,
% by the properties of Wasserstein distance and
from Theorem~\ref{mainthm}, if $n\bR> r_0$,\vspace*{1pt}
\[
d_W\bigl(\cL(W),\NB(n\bR,\th_T)\bigr) \Le
\frac{{24}(1+\log\{1/(1-\th_T)\}
)\sqrt{\th_T}}{(1-\theta_T)^{3/2}\sqrt{n\bR}} \sum_{i=1}^n
\bigl(A_T^*\bigr)^{(i)},
\]
where $r_0$ is as for \Ref{ADB-r_0-bnd}.
% \frac{3}{(1-\th_T)^4\sqrt{2R_a}}
Defining $\s:= n^{-1}\sum_{i=1}^n(A_T^*)^{(i)}$,
the bound grows with $n$ roughly as $\s\sqrt n/\bR$. However, the variability
of the distribution $\NB(n\bR,\th_T)$ is also on the scale $\sqrt
n$, so that the relevant measure
of distance is $n^{-1/2}d_W(\cL(W),\NB(n\bR,\th_T))$, which is
small provided that $\s\ll\bR$.
If $\NB(n\bR,\th_T)$ is replaced by $\NB(nR_a^*,\th_T)$, the
additional term $|\sum_{i=1}^nA_T^{(i)}|$
in $d_W(\cL(W),\NB(nR_a^*,\th_T))$ is
also roughly of order $\s\sqrt n$, if, for instance, the $A_T^{(i)}$
are independent random variables
with mean zero.\vspace*{1pt}

%sA #&#
\begin{appendix}\label{appendix}
\section*{Appendix}
\setcounter{equation}{10}
\setcounter{section}{3}
\renewcommand{\theequation}{\arabic{section}.\arabic{equation}}

The constant $K$ in \Ref{ADB-bound-with-delta} can be shown to
satisfy $K\le16$ as follows.
Expression \Ref{f_j-dash} can be written in a neat form:\vspace*{1pt}
\begin{eqnarray*}
f_j'(\th) &\Eq&\th^{j-3}(1-\th)\bigl\{j(j-1)
\bigl((\th_T-\th)^2+(\th_T-\th ) (1-
\th_T)\bigr)
\\[2pt]
&&\phantom{\th^{j-3}(1-\th)\bigl\{}{}-2j\bigl(\th_T-\th^2\bigr)+2
\th_T\bigr\},
\end{eqnarray*}
from which it follows that\vspace*{1pt}
%
%
%eA.11 #&#
\begin{equation}
\label{xiaadd2} \bigl |f_2'(\th)\bigr | \Le2(1-\th) (3\th+1+
\th_T)
\end{equation}
and, for $j\ge3$ and $0 \le\th\le\th_T$,\vspace*{1pt}
%
%
%eA.12 #&#
\begin{eqnarray}
\label{xiaadd3} \bigl |f_j'(\th)\bigr | &\Le&\th^{j-3}(1-\th)
\bigl\{j(j-1) \bigl((\th_T-\th)^2+(\th_T-\th )
(1-\th_T)\bigr)
\nonumber
\\[-7pt]
\\[-7pt]
&&\phantom{\th^{j-3}(1-\th)
\bigl\{}{}+2j\bigl(\th_T-\th^2\bigr)+2
\th_T\bigr\}.  \nonumber
\end{eqnarray}
Now, \Ref{xiaadd2} yields\vspace*{1pt}
\[
\int_0^{\th_T}\bigl |f_2'(\th)\bigr |
\,\mathrm{d}\th\le \int_0^{\th
_T}2(1-\th) (3\th +1+
\th_T)\,\mathrm{d}\th\Eq % &=&
4\th_T^2+2
\th_T-3\th_T^3,
\]
and, for $j\ge3$, integrating \Ref{xiaadd3} gives\vspace*{1pt}
\[
\int_0^{\th_T}\bigl |f_j'(\th)\bigr |
\,\mathrm{d}\th % &&\le\int_0^{\th_T}\th^{j-3}(1-\th)\{j(j-1)((\th_T-\th)^2+(\th_T-
\Le3\th_T^{j-1}(1-
\th_T)^2+4\th_T^{j-1} \biggl(
\frac{2}{j-2}-\frac
{\th_T^2}{j+1}-\frac{\th_T}{j-1} \biggr).
\]
Hence,
%
%
%eA.13 #&#
\begin{eqnarray}
\label{xiaadd4} && \sum_{j\ge2}(j-1)\int
_0^{\th_T}\bigl |f_j'(\th)\bigr | \,
\mathrm{d}\th
\nonumber
\\
&&\quad\Le\sum_{j\ge3}(j-1) \biggl\{3
\th_T^{j-1}(1-\th_T)^2 +4
\th_T^{j-1} \biggl(\frac{2}{j-2}-\frac{\th_T^2}{j+1}-
\frac{\th
_T}{j-1} \biggr) \biggr\}
\nonumber
\\
&& \qquad{} +4\th_T^2+2\th_T-3
\th_T^3
\\
&&\quad\Eq-6\th_T+14\th_T^2-(14/3)
\th_T^3-8(\th_T+1)\ln(1-\th_T)
\nonumber
\\
&&\quad\Le2\th_T+14\th_T^2-(14/3)
\th_T^3-16\th_T\ln(1-\th _T)
\nonumber
\end{eqnarray}
and
%
%
%eA.14 #&#
\begin{equation}
\label{xiaadd5} \sum_{j\ge2}(j-1)\bigl |A_Tf_j(
\th_T)\bigr | \Eq|A_T|\sum_{j\ge2}(j-1)
(1-\th _T)^2\th_T^{j-1}
\Eq|A_T|\th_T.
\end{equation}
Combining \Ref{ADB-pause}, \Ref{ADB-integrated-j-bnd}, \Ref
{ADB-second-bnd}, \Ref{xiaadd4} and \Ref{xiaadd5} yields
\begin{eqnarray*}
&&\bigl\llvert \E \bigl(\theta_T(R_T+W)g(W+1)-Wg(W)
\bigr)\bigr\rrvert
\\
&&\quad\Le\|\Delta g\| \bigl\{3\th_T+14\th_T^2-(14/3)
\th _T^3-16\th _T\ln(1-\th_T)
\bigr\}
\\
&&\quad\Le\|\Delta g\| \th_T A_T^*\bigl(37/3 +16 \ln
\bigl\{1/(1-\th_T)\bigr\}\bigr).
\end{eqnarray*}
\end{appendix}

\section*{Acknowledgements}

Work supported in part by Australian Research Council Grants Nos
DP120102728 and DP120102398.

%suskaldyti doi

%
% imsref loaded by audrone.aklyte, 2014-03-20 12:18:27
%

\printhistory

\end{document}